\documentclass[review]{elsarticle}
\usepackage{lineno,hyperref}
\modulolinenumbers[5]

\usepackage{array}
\usepackage{graphicx}      
\usepackage{caption}        
\usepackage{algorithm}
\usepackage{algorithmic}
\usepackage{subcaption}     
\usepackage{amsmath,amsfonts,amssymb}
\usepackage{epsfig,makecell,float}
\usepackage{setspace,mathrsfs}
\usepackage{tocloft}
\usepackage{textcomp}
\usepackage{multirow,indentfirst,times,color}
\usepackage{booktabs}
\usepackage{threeparttable}
\usepackage{amsthm}
\usepackage{extarrows}

\usepackage[titletoc]{appendix}
\usepackage{epstopdf}

\usepackage{natbib}
\biboptions{numbers,sort&compress} 
\usepackage{bm}

\makeatletter
\newif\if@restonecol
\makeatother

	\newtheorem{thm}{Theorem}

\let\clearpage\relax

\journal{Journal of Digital Signal Processing}

\begin{document}\textbf{}

\begin{frontmatter}

\title{A novel STAP algorithm via volume cross-correlation function on the Grassmann manifold}

\author[mymainaddress,mysecondaryaddress]{Jia-Mian Li}

\author[mymainaddress,mysecondaryaddress]{Jian-Yi Chen}

\author[mymainaddress,mysecondaryaddress]{Bing-Zhao Li\corref{mycorrespondingauthor}}
\cortext[mycorrespondingauthor]{Corresponding author}
\ead{li_bingzhao@bit.edu.cn}

\address[mymainaddress]{School of Mathematics and Statistics, Beijing Institute of Technology, Beijing 100081, China}
\address[mysecondaryaddress]{Beijing Key Laboratory on MCAACI, Beijing Institute of Technology, Beijing 100081, China}

\begin{abstract}
The performance of space-time adaptive processing (STAP) is often degraded by factors such as limited sample size and moving targets. Traditional clutter covariance matrix (CCM) estimation relies on Euclidean metrics, which fail to capture the intrinsic geometric and structural properties of the covariance matrix, thus limiting the utilization of structural information in the data.
To address these issues, the proposed algorithm begins by constructing Toeplitz Hermitian positive definite (THPD) matrices from the training samples. The Brauer disc (BD) theorem is then employed to filter out THPD matrices containing target signals, retaining only clutter-related matrices. These clutter matrices undergo eigendecomposition to construct the Grassmann manifold, enabling CCM estimation through the volume cross-correlation function (VCF) and gradient descent method. 
Finally, the filter weight vector is computed for filtering. By fully leveraging the structural information in radar data, this approach significantly enhances both accuracy and robustness of clutter suppression. 
Experimental results on simulated and measured data demonstrate superior performance of the proposed algorithm in heterogeneous environments.
\end{abstract}



\end{frontmatter}


\section{Introduction}
Clutter significantly affects the target detection and identification capabilities of radar systems \cite{RN68,RN50,RN65}. Therefore, effectively suppressing sea clutter is crucial for improving the accuracy and reliability of radar systems \cite{RN73,RN75}. 
Space-time adaptive processing (STAP) leverages the spatial and temporal dimensions of radar echoes, allowing for adaptive adjustment of filter parameters based on different environments and mission requirements to suppress clutter \cite{RN57,RN79,RN81}.
Traditional STAP methods often perform poorly in non-homogeneous clutter environments because they typically use a sample covariance matrix (SCM) as a substitute for the clutter covariance matrix (CCM). When the sample size is insufficient, the estimation of the SCM becomes inaccurate, affecting STAP performance \cite{RN84,RN59,RN54}. In the presence of range-spread targets, which occupy multiple range cells, training samples may contain multiple target signals, leading to inaccurate estimation of the CCM and target self-nulling, which reduces suppression effectiveness \cite{RN62,RN53,RN61}. Therefore, achieving more accurate CCM estimation in the presence of limited samples and target interference is a key challenge.

Current methods include using the Generalized Inner Product (GIP) technique \cite{RN95} to filter out training samples contaminated by target signals, and employing Random Matrix Theory (RMT) with the spiked covariance model to estimate the matrix by determining the eigenvalues of the SCM \cite{RN77}. The knowledge-aided STAP method utilizes prior knowledge to improve non-homogeneous clutter suppression performance. This method guides the selection of training strategies and samples or constructs the CCM, which is then fused with the estimated CCM to form the final covariance matrix \cite{RN66,RN64,RN51,RN72}. 
Reference \cite{RN86} introduces a knowledge-aided Bayesian framework in which covariance matrices are considered as random variables, with partial information about the covariance matrix of the training samples available. In this context, the maximum a posteriori (MAP) estimate for the primary data covariance matrix is derived.
The sparse recovery STAP method \cite{RN71} has become a research hotspot in radar technology, improving clutter suppression performance for airborne radar under conditions of severe sample insufficiency. However, this method also faces performance degradation due to lack of detailed prior knowledge. To address this, \cite{RN94} proposes an STAP algorithm based on structured two-level block sparsity, which improves clutter suppression and target detection performance with limited training samples. Another study proposes a STAP filter based on sparse Bayesian learning theory \cite{RN48}, achieving excellent clutter suppression and target detection performance. Furthermore, \cite{RN74} introduces a robust STAP filter based on the volume cross-correlation function (VCF), which improves clutter suppression under unknown clutter distributions. This filter does not require prior information about the clutter distribution, demonstrating robustness in cluttered environments.

By transforming radar data into Toeplitz Hermitian Positive Definite (THPD) matrices \cite{RN92,RN93}, smoothing estimations can be achieved to reduce variance and improve stability, which helps enhance detection and estimation performance. However, these methods have limitations \cite{RN46}: traditional CCM estimation is based on Euclidean metrics, which cannot always maintain positive definiteness under addition and scalar multiplication, making them unsuitable for processing THPD matrices. Moreover, these methods overlook the intrinsic geometric properties and structural features of the covariance matrix, failing to fully exploit the structural information of the data.

Therefore, to address the shortcomings of existing technologies and achieve better clutter suppression, we propose a novel algorithm based on Brauer disc (BD) theorem \cite{RN98,RN99} and VCF on the Grassmann manifold (GVCF), referred to as BGVCF.
The main contributions are as follows. First, sliding windows are used to determine the test cells, guard cells, and training samples, and the THPD matrix of the training samples is calculated. Next, based on the BD theorem, clutter THPD matrices without target contamination are selected. The matrices are then eigendecomposed to form the Grassmann manifold\cite{RN104,RN101}, and the VCF, along with gradient descent, is applied to estimate the CCM. Finally, a space-time steering vector is constructed, and the filter weight vector is calculated. Compared with existing methods, this approach can more accurately and effectively distinguish target signals from clutter, ensuring that covariance matrix estimation is not affected by target interference. Moreover, the method fully utilizes the intrinsic geometric properties and structural features of the covariance matrix, does not rely on specific clutter distribution assumptions, and achieves robust covariance matrix estimation by directly measuring the distance between the sampled signal subspace and the target subspace, leading to effective clutter suppression and significantly improving suppression accuracy and stability.

The remainder of this paper is organized as follows. Section \ref{sec2} introduces the STAP signal model and related theory. Section \ref{sec3} provides a detailed implementation of the BGVCF-STAP algorithm. 
Section \ref{sec4} reports numerical experiments using both simulated and measured data, demonstrating the effectiveness of the proposed method in clutter suppression. Finally, Section \ref{sec5} offers concluding remarks and discussions.

\section{Signal model and related theory}
\label{sec2}
\subsection{Signal model}
In this section, we examine an airborne pulse-Doppler radar system equipped with a uniform linear array (ULA) consisting of \(M\) elements, each spaced equally by a distance \(d\). During a coherent processing interval (CPI), each element receives \(N\) pulses at a constant pulse repetition frequency (PRF). Target detection is framed as a binary hypothesis testing problem \cite{RN82}. Under hypothesis \(H_1\), the received signal includes both the target echo and clutter noise, while under hypothesis \(H_0\), only clutter and noise are present. The space-time snapshot data \( \mathbf{x}(l) \in \mathbb{C}^{MN \times 1} \) collected from the \(l\)th range cell is represented as
\begin{align}
\begin{cases}
H_0: \mathbf{x}(l) = \mathbf{c}(l) + \mathbf{n}(l) \\
H_1: \mathbf{x}(l) = \mathbf{s}(l) + \mathbf{c}(l) + \mathbf{n}(l)
\end{cases}
\end{align}
where \( \mathbf{s}(l) \) is the target space-time steering vector, \( \mathbf{c}(l) \) is the clutter space-time steering vector, \( \mathbf{n}(l) \) is complex Gaussian noise with a variance of \( \delta^2 \).
Meanwhile, the clutter signal \( \mathbf{c}(l) \) can be expressed as \cite{RN81}
\begin{align}
\mathbf{c}(l) = \sum_{j=1}^{N_{a}} \sum_{i=1}^{N_{c}} a_{c,i,j} \boldsymbol{v}\left(f_{d,c,i,j}, f_{s,c,i,j}\right)
\end{align}
where \( N_a \) are the number of range ambiguities, \( N_c \) are the numbers of clutter patches within each range blur, the coefficient \( a_{c,i,j} \) represents the complex scattering amplitude of the clutter patch, \( f_{d,c,i,j} \) and \( f_{s,c,i,j} \) are the corresponding normalized Doppler and spatial frequencies.
$\boldsymbol{v}\left(f_{d,c,i,j}, f_{s,c,i,j}\right)=\boldsymbol{v}_{d}\left(f_{d,c,i}\right) \otimes \boldsymbol{v}_{s}\left(f_{s,c,i}\right)$ is the space-time oriented vector of the $i$th clutter scattering point, $\otimes$ representing the Kronecker product, $\boldsymbol{v}_d\left(f_{d,c,i}\right)$ denotes the time-steering vector and $\boldsymbol{v}_s\left(f_{s,c,i}\right)$ denotes the space-steering vector, which are defined as
\begin{equation}
\begin{aligned}
\label{(3)}
\bm{v}_d\left(f_{d,c,i}\right) &= \left[1, \exp\left(j 2\pi f_{d,c,i}\right), \dots, \exp\left(j 2\pi (M-1) f_{d,c,i}\right)\right]^T \\
\bm{v}_s\left(f_{s,c,i}\right) &= \left[1, \exp\left(j 2\pi f_{s,c,i}\right), \dots, \exp\left(j 2\pi (N-1) f_{s,c,i}\right)\right]^T
\end{aligned}
\end{equation}

The snapshot of target \( \mathbf{s}(l) \) can be denoted as
\begin{align}
\mathbf{s}(l) = a_{t}\boldsymbol{v}\left(f_{d,t}, f_{s,t}\right)
\end{align}
where \( a_{t} \) is the reflection coefficient of the target, $\boldsymbol{v}\left(f_{d,{t}}, f_{s,{t}}\right)=\boldsymbol{v}_{d}\left(f_{d, {ti}}\right) \otimes \boldsymbol{v}_{s}\left(f_{s, {ti}}\right)$ is the space-time steering vector of the target, its expression is the same as (\ref{(3)}).

The goal of the space-time filter is to minimize the filtered output power of interference 
while ensuring that the output power at the target position is not lost. Specifically, 
the design of the weight vector for the space-time filter can be obtained by solving the following optimization problem
\begin{align}
\min_{\mathbf{w}} \mathbf{w}^H \mathbf{R} \mathbf{w}
\quad \text{s.t.} \quad \mathbf{w}^H \bm{v}(f_{d,t}, f_{s,t}) = 1
\end{align}
where the weight vector of the STAP filter is given by
\begin{align}
\label{(6)}
\mathbf{w}_{opt} = \frac{\mathbf{R}^{-1} \bm{v}(f_{d,t}, f_{s,t})}
{\bm{v}(f_{d,t}, f_{s,t})^H \mathbf{R}^{-1} \bm{v}(f_{d,t}, f_{s,t})}
\end{align}
where  
\( \mathbf{R} \) is the ideal CCM and \( [\cdot]^H \) denotes the conjugate transpose operation. In practical applications, the ideal CCM is usually unknown and needs to be estimated through statistical methods. The practical approximation of the CCM can be achieved by averaging the autocorrelation of non-target range cell data, as shown below
\begin{align}
\hat{\mathbf{R}} = \frac{1}{L} \sum_{l=1}^{L} \mathbf{x}(l) \mathbf{x}^H(l)
\end{align}
where \( L \) is the number of training samples. However, in heterogeneous environments, this estimation may not be accurate.

Typically, \( \hat{\mathbf{R}} \) is used in place of \( \mathbf{R} \) in (\ref{(6)}) to calculate the adaptive filtering weight \( \mathbf{w} \). 
The weighted processing is carried out by taking the inner product of the space-time snapshot data from the current range cell with the corresponding weight vector, which is given by
\begin{align}
y_l = \mathbf{w}^H \mathbf{x}(l)
\end{align}

\subsection{Toeplitz Hermitian positive definite matrix}

The concept of a Toeplitz matrix is first introduced. Consider \( T = \{t_{ij}\} \in \mathbb{C}^{n \times n} \), where \( t_{ij} = t_{i-j} \) for \( i, j = 1, 2, \ldots, n \). A matrix \( T \) satisfying this condition is referred to as a Toeplitz matrix. Following this, the THPD covariance matrix is defined as a matrix that combines three properties: it is a Toeplitz matrix, a Hermitian matrix, and positive definite. 
The THPD covariance matrix C has the following \cite{barbar}
\begin{align}
C=\begin{bmatrix}c_0&c_{-1}&c_{-2}&\cdots&c_{-n+1}\\c_1&c_0&c_{-1}&\cdots&c_{-n+2}\\c_2&c_1&c_0&\cdots&c_{-n+3}\\\vdots&\vdots&\vdots&\ddots&\vdots\\c_{n-1}&c_{n-2}&c_{n-3}&\cdots&c_0\end{bmatrix}
\end{align}
where $c_i = \overline{c_{-i}}$ to guarantee the Hermitian property, and the matrix is positive definite. 

The matrix is derived from the calculation of the reflection coefficient. $\mathbf{x}$ represent a snapshot containing \( N \) elements. Assuming the signal follows an autoregressive (AR) model, the reflection coefficient corresponds to the \( k \)th order regression coefficient. This coefficient can be computed using the regularized Burg algorithm, as detailed below \cite{RN92}
\begin{align}
    f_0(k) &= b_0(k) = \mathbf{x}(k), \quad k = 0, \ldots, N - 1 \\
    \label{(11)}
    P_0 &= \frac{1}{N} \sum_{k=0}^{N-1} |\mathbf{x}(k)|^2
\end{align}
for $n=1$ to $N-1$ do
\begin{align}
\label{(12)}
\mu_{n} &= 
-\frac{
    \frac{2}{N - n} \sum\limits_{k = n + 1}^{N} f_{n-1}(k) \bar{b}_{n-1}(k - 1) 
    + 2 \sum\limits_{k=1}^{n-1} \nu_k^{(n)} a_{k}^{(n-1)} a_{n-k}^{(n-1)}
}{
    \frac{1}{N - n} \sum\limits_{k = n + 1}^{N} \left( |f_{n-1}(k)|^{2} + |b_{n-1}(k-1)|^{2} \right) 
    + 2 \sum\limits_{k=0}^{n-1} \nu_k^{(n)} |a_{k}^{(n-1)}|^2
}
\end{align}
where
\begin{align}
&\nu_{k}^{(n)} =\psi_{1}(2\pi)^{2}(k-n)^{2}\\
&\begin{cases}a_0^{(n)}=1\\a_k^{(n)}=a_k^{(n-1)}+\mu_n\tilde{a}_{n-k}^{(n-1)},\quad k=1,\ldots,n-1\\a_n^{(n)}=\mu_n\end{cases}
\end{align}
and
\begin{align}
&\begin{cases}f_n(k)=f_{n-1}(k)+\mu_nb_{n-1}(k-1)\\b_n(k)=b_{n-1}(k-1)+\bar{\mu}_nf_{n-1}(k)\end{cases}
\end{align}

\subsection{VCF on the Grassmann manifold}
The union model of linear subspaces is a commonly used signal model. 
The signal $\mathbf{x}$ in this model is assumed to be a vector in the union of linear subspaces in $\mathbb{R}^N$, defined as
\begin{align}
\label{(16)}
\mathbf{x} \in \bigcup_{i=1}^{L} X_i, \quad X_i = \left\{ \mathbf{x} = X_i \mathbf{a}_i \mid X_i \in \mathbb{R}^{N \times k}, \mathbf{a}_i \in \mathbb{R}^k \right\}
\end{align}
where the column vectors of the matrix $X_i$ form the basis of the corresponding subspace $X_i$, and $\dim(X_i) = k < N$. 

The Grassmann manifold $Gr(k, N)$ is a topological space \cite{RN104,RN101} in which each point corresponds to a $k$-dimensional linear subspace of $\mathbb{R}^N$. 
In general, the union of subspaces in (\ref{(16)}) is equivalent to a finite set of points in $Gr(k, N)$, namely
\begin{align}
\label{(17)}
\zeta(k, N, L) = \{ X_1, \dots, X_L \}
\end{align}
where $X_i \in Gr(k, N), \, i = 1, \dots, L$ represents a point on the Grassmann manifold $Gr(k, N)$. 

Volume, as a significant parameter in the Grassmann manifold space, quantifies the degree of separation between two linear subspaces and is closely related to the principal angles between them. Based on the geometric structure of the subspaces, the volume is utilized to derive a VCF, which characterizes the distance between subspaces.
The $d$-dimensional volume of a full rank matrix $S$ is defined as \cite{RN109}
\begin{align}
\mathrm{vol}_{d}\big(S\big)=\prod_{i=1}^{d}\sigma_i
\end{align}
where $\sigma_i, i = 1, \dots, d$ represent the singular values of the matrix $S$.
Further, for two subspaces whose basis matrices are denoted by $S_{1}\in\mathbb{C}^{{K}\times d_{1}}$, $S_{2}\in\mathbb{C}^{{K}\times d_{2}}$, then their volume-dependent functions are defined as
\begin{align}
\begin{gathered}
\mathrm{VCF}\left(S_{1},S_{2}\right)=\frac{\mathrm{vol}_{d_{1}+d_{2}}\left(\left[S_{1},S_{2}\right]\right)}{\mathrm{vol}_{d_{1}}\left(S_{1}\right)\mathrm{vol}_{d_{2}}\left(S_{2}\right)} \\
=\prod_{i=1}^{\min(d_1,d_2)}\sin\theta_i\left(span(S_1),span(S_2)\right) 
\end{gathered}
\end{align}
where $\theta_1, \theta_2, \dots, \theta_{\min(d_1, d_2)}$ are the principal angles between the two subspaces spanned by $S_1$ and $S_2$, respectively. 
The physical meaning is the angle between the $d_1$ orthonormal vectors in the subspace $span(S_1)$ and the $d_2$ orthonormal vectors in the subspace $span(S_2)$. 

\section{Proposed BGVCF-STAP algorithm}
\label{sec3}
\begin{figure}
	\centering
	\includegraphics[width=1\textwidth]{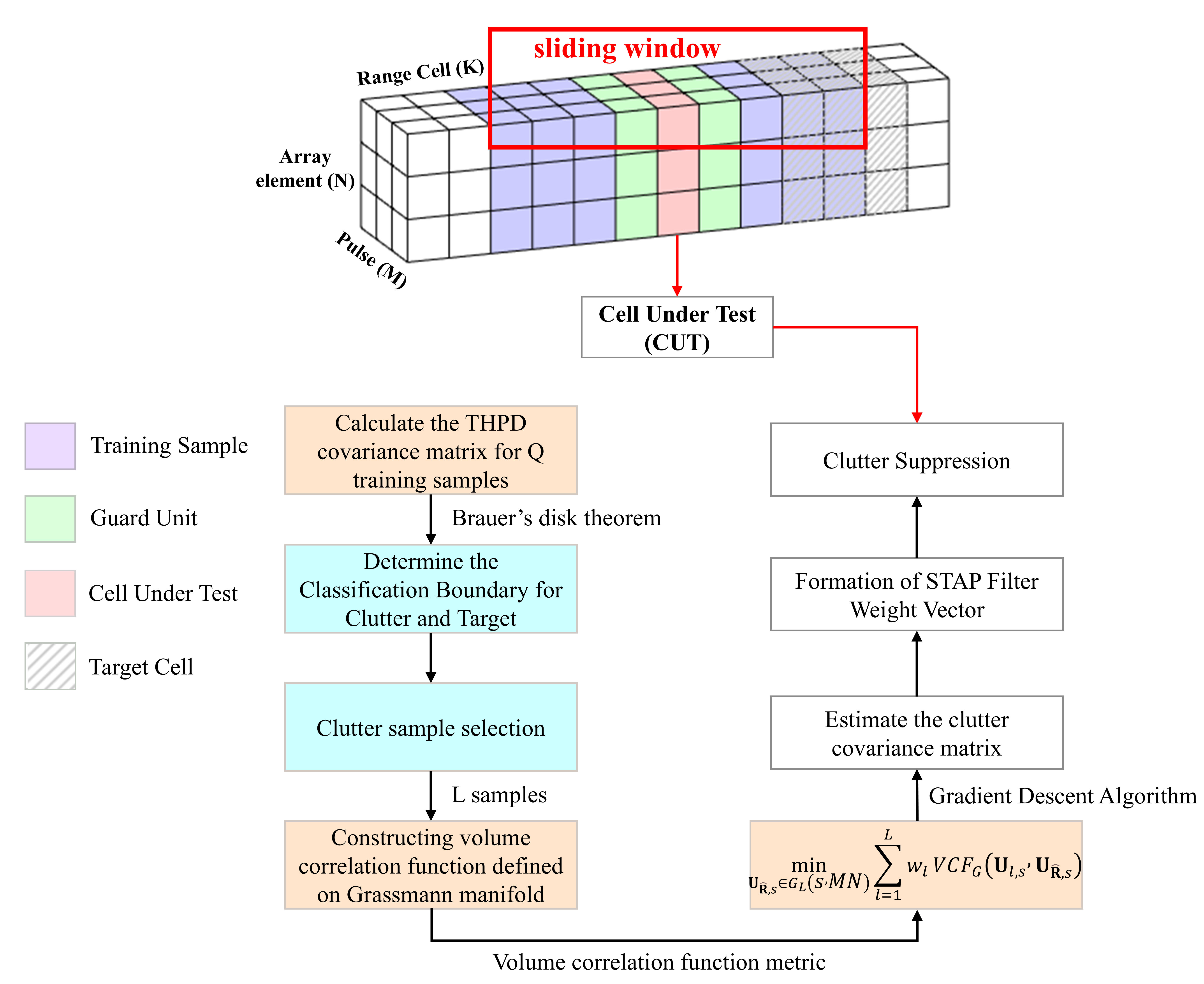}
	\caption{Algorithm flowchart. }
 \label{fig1}
\end{figure}

This section provides a detailed description of the BGVCF-STAP algorithm, with the flowchart shown in Fig. ~\ref{fig1}. The BGVCF algorithm is primarily divided into three stages: the calculation of the THPD matrix for training samples, the clutter and target screening stage, and the CCM estimation stage.
\subsection{Calculation of the THPD matrix}
For each range cell, a sliding window of fixed length is applied along the range dimension to form the training sample matrix from the received radar data.  
The $q$th training sample is given by
$x(q) = \left[ x_q(1) \; x_q(2) \; \cdots \; x_q(MN) \right]^T,$
where $q = 1, 2, \cdots, Q$. The training sample matrix is selected using a sliding window of length $Q$, excluding the test cell and guard cells.

Assume that the covariance matrix of clutter has the following THPD structure, which is denoted as the THPD covariance matrix, namely.  
\begin{align}
\mathbf{R}_q = 
\begin{bmatrix}
    r_0 & r_1^* & \cdots & r_{MN-1}^* \\
    r_1 & r_0 & \cdots & r_{MN-2}^* \\
    \vdots & \vdots & \ddots & \vdots \\
    r_{MN-2} & r_{MN-3} & \cdots & r_1^* \\
    r_{MN-1} & r_{MN-2} & \cdots & r_0
\end{bmatrix}
\end{align}
where each $r_i$ can be calculated using the reflection coefficient $\mu_i$, as given by
\begin{align}
\begin{cases}r_0&=P_0,\quad r_1=-P_0\mu_1\\r_i&=-P_{i-1}\mu_i\boldsymbol{\alpha}_{i-1}^T\mathbf{J}_{i-1}\mathbf{R}_{i-1}^{-1}\boldsymbol{\alpha}_{i-1},\quad2\leqslant i\leqslant MN-1\end{cases}
\end{align}
$P_0$ and $\mu_i$ are obtained through the regularized Burg algorithm (\ref{(12)})–(\ref{(17)}).
\begin{align}
&\alpha_{i-1} =
\begin{bmatrix}
r_1 \\
\vdots \\
r_{i-1}
\end{bmatrix},
\quad
\mathbf{J}_0 = 1, \quad
\mathbf{J}_1 =
\begin{bmatrix}
0 & 1 \\
1 & 0
\end{bmatrix}, \\
&\mathbf{J}_{i-1} =
\begin{bmatrix}
0 & 0 & \cdots & 1 \\
0 & \cdots & 1 & 0 \\
\vdots & \cdots & \cdots & \vdots \\
1 & \cdots & 0 & 0
\end{bmatrix}, \\
&\mathbf{R}_{i-1} = \mathbf{R}(2:i, 1:k-1), \quad
P_{i-1} = P_0 \prod_{j=1}^{i-1} \left(1 - |\mu_j|^2\right).
\end{align}
Through the above definition, the THPD covariance matrix of Q training samples is calculated.

\subsection{Distinguish clutter-plus-noise and targets}
\label{3.2}
\begin{thm}
Brauer Disk Theorem \cite{RN98}: Suppose $R = [r_{ij}]$ is an $n$-th order complex Toeplitz matrix, with $r_{ii} = r_0$, and $n \geq 2$. Let $\sigma(R)$ denote the eigenvalues of $R$. Then,
\begin{align}
\sigma\left(R\right) &\subseteq \mathrm{Y}\left(R\right) = \bigcup_{i=1}^{\left\lceil\frac{n}{2}\right\rceil} \mathrm{Y}_{i}\left(R\right) \\
\mathrm{Y}_{i}\left(R\right) &= \left\{z \in C : \left|z - r_{0}\right| \leq D\right\} \\
D &= \operatorname*{max}_{i,j \in N, i \neq j} \sqrt{a_{i}\left(R\right) a_{j}\left(R\right)} \\
a_{i}\left(R\right) &= \sum_{i \neq j} \left|r_{ij}\right| \\
\left\lceil \frac{n}{2} \right\rceil &=
\begin{cases}
\frac{n}{2}, & n = 2k, \, k \in \mathbb{Z} \\
\frac{n+1}{2}, & n = 2k+1, \, k \in \mathbb{Z}
\end{cases}
\end{align}
where $\mathrm{Y}\left(R\right)$ is the union of the Brauer circle sets of the symmetric matrix $R$.
\end{thm}

The BD theorem is a mathematical tool \cite{RN99} used to determine the inclusion region of eigenvalues for a square matrix. It takes into account both the sum of the elements in the rows and the sum of the elements in the columns of the matrix. This enables it to define the location of the eigenvalues more precisely, providing tighter bounds.
Based on BD theorem, the clutter and target signal THPD covariance matrix are classified. The basic principle is as follows:

\begin{enumerate}[(1)]
\item \text{Eigenvalue analysis:}
each THPD covariance matrix is decomposed by eigenvalue to obtain its eigenvalue. The center of the largest BD, the contained region of the largest eigenvalues, $P_0$ represents the snapshot power, as shown in (\ref{(11)}).

\item \text{Determination of the distribution range of eigenvalues:}
according to the distribution range of eigenvalues, the minimum BD radius $\rho$  of the maximum eigenvalue of each THPD covariance matrix is calculated. This radius describes the distribution range of eigenvalues in the complex plane, reflects the degree of disturbance of the matrix, and represents the minimum limit of clutter-plus-noise.

\item \text{Determination of classification limits:}
using the BD centers of all THPD covariance matrices, the Brauer cluster boundary limit $T_B$ for clutter-plus-noise can be established \cite{RN45}
\begin{align}
T_B=\frac{\frac{1}{M}\sum_{k=1}^K\hat{\mathbf{R}}_k}{\left(\prod_{k=1}^K\hat{\mathbf{R}}_k\right)^{\frac{1}{M}}}\rho 
\end{align}
where $\frac{1}{M} \sum_{k=1}^K \hat{\mathbf{R}}_k \; \text{and} \: 
\left( \prod_{k=1}^K \hat{\mathbf{R}}_k \right)^{\frac{1}{M}}$
represent the arithmetic mean and geometric mean of the BD centers, respectively. In general, the eigenvalue distribution range of a clutter matrix is relatively small, resulting in a smaller minimum BD radius. However, for matrices containing target signals, the situation is the opposite.

\item \text{Classification and screening:}
The minimum BD radius of all THPD covariance matrices is compared with the classification threshold $T_B$, and matrices with a radius smaller than $T_B$ are considered as clutter matrices. Specifically, the THPD matrices corresponding to eigenvalues within the Brauer cluster boundary are regarded as clutter-plus-noise covariance matrices, while those corresponding to eigenvalues outside the Brauer cluster boundary are considered as target matrices.
\end{enumerate}

\subsection{Estimation of CCM based on GVCF}
In general, the CCM estimation for clutter suppression can be viewed as a function of the matrix composed of the received data in the training data set
\begin{align}
\hat{\mathbf{R}}=\phi(\hat{\mathbf{R}}_1,\hat{\mathbf{R}}_q,\cdots,\hat{\mathbf{R}}_Q)=\phi'(\boldsymbol{x}_1,\boldsymbol{x}_q,\cdots,\boldsymbol{x}_Q)
\end{align}
where $\hat{\mathbf{R}}_q$ represent the THPD covariance matrix of the $q$th training sample, and $\hat{\mathbf{R}}$ be the estimated CCM. $x_q = x(q)$ denotes the training data without the target. The CCM can be estimated using the following optimization model
\begin{align}
\min_{\hat{\mathbf{R}} \in \mathbb{R}} \sum_{q=1}^{Q} w_q d_E^2(\hat{\mathbf{R}}_q, \hat{\mathbf{R}})
\end{align}
where $\mathbb{R}$ represents the linear space spanned by the observation matrices, $w_q > 0$ are the weights that satisfy $\sum_{q=1}^{Q} w_q = 1$, indicating the proportion of each training unit in the estimation. $d_E^2(\hat{\mathbf{R}}_q, \hat{\mathbf{R}})$ represents the Euclidean distance between the two spaces.
\begin{figure}
	\centering
	\includegraphics[width=0.5\textwidth]{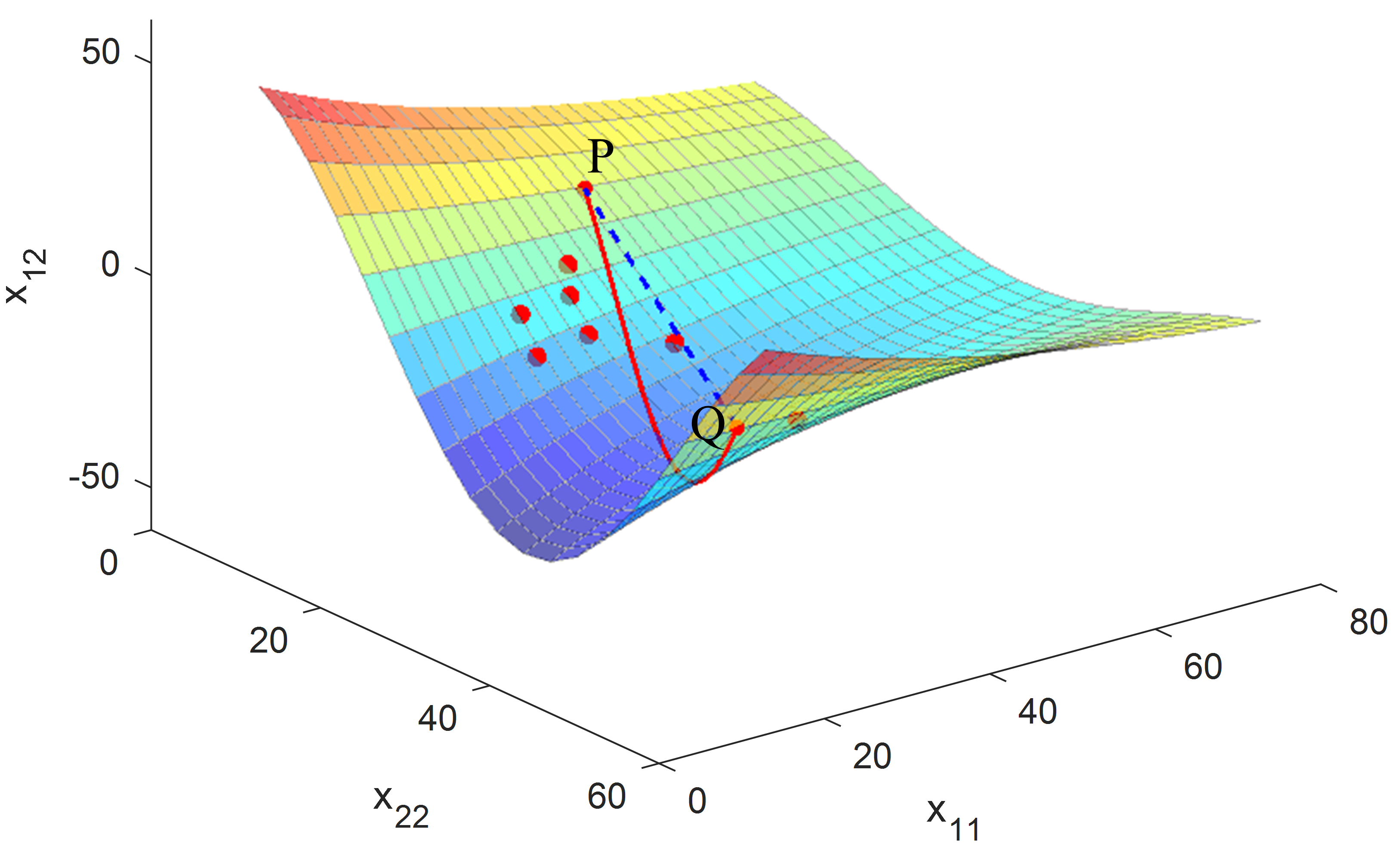}
	\caption{Schematic diagram of distance measurement on manifolds. }
 \label{fig2}
\end{figure}

However, the Euclidean distance only focuses on the signal energy and does not fully exploit the structural information in the data. Specifically, the Euclidean distance $d_E(\hat{\mathbf{R}}_q, \hat{\mathbf{R}})$ calculates the discrepancy between two matrices using the Frobenius norm, which is defined as
\begin{align}
d_E^2(\hat{\mathbf{R}}_q,\hat{\mathbf{R}})=||\hat{\mathbf{R}}_q-\hat{\mathbf{R}}||_F^2=\mathrm{Tr}\left((\hat{\mathbf{R}}_q-\hat{\mathbf{R}})(\hat{\mathbf{R}}_q-\hat{\mathbf{R}})^H\right)
\end{align}
where \( \| \mathbf{\cdot} \|_F \) denotes the Frobenius norm, \( \mathrm{Tr}(\mathbf{\cdot}) \) represents the trace of the matrix, 
this measure overlooks the internal geometric and structural properties of the CCM.

Fig. ~\ref{fig2} illustrates distance measurement on a manifold, highlighting the limitations of the Euclidean metric in high-dimensional structures due to its neglect of geometric information \cite{RN102}. The red circles represent positive definite covariance matrices, with the axes \( x_{11} \), \( x_{22} \), and \( x_{12} \) visualizing the matrix space. In the figure, the Euclidean distance (blue dashed line) fails to accurately reflect the true distance between points P and Q on the manifold, while the red geodesic path represents the actual intrinsic distance \cite{RN103}.
The limitation of the Euclidean distance lies in its inability to preserve the positive definiteness of matrices. For instance, when a positive definite matrix is scaled by a negative scalar, the result is no longer positive definite. Therefore, using the Euclidean distance to measure THPD covariance matrices is inappropriate, whereas the structure of the Grassmann manifold naturally preserves this characteristic.

The Grassmann manifold $G(k, n)$ is the collection of all $k$-dimensional subspaces within an $n$-dimensional space. In other words, the Grassmann manifold is the set of all $n \times k$ matrices, where each $n \times k$ matrix represents a $k$-dimensional subspace, and these matrices are orthonormal. For $Q$ THPD matrices of size $MN \times MN$, eigenvalue decomposition or singular value decomposition is performed. Since $\hat{\mathbf{R}}_q$ is Hermitian positive definite, it can be expressed as
\begin{align}
\hat{\mathbf{R}}_q = \mathbf{U}_q \bm{\Lambda}_q \mathbf{U}_q^H
\end{align}
where $\mathbf{U}_q$ is the eigenvector matrix of the $q$-th THPD matrix of size $MN \times MN$, and $\bm{\Lambda}_q$ is the diagonal matrix of eigenvalues of the $q$-th THPD matrix. Each diagonal element corresponds to an eigenvalue. From each $\mathbf{U}_q$, the first $s$ eigenvectors (corresponding to the largest $s$ eigenvalues) are selected to form an $MN \times s$ matrix $\mathbf{U}_{q, s}$. The column space of this matrix is an $s$-dimensional subspace.
Each THPD matrix $\hat{\mathbf{R}}_q$ is mapped to an $s$-dimensional subspace, and the collection of these subspaces forms a point on the Grassmann manifold $G_Q(s, MN)$.

Based on the structure of the Grassmann manifold $G_Q(s, MN)$, the volume cross-correlation function on the Grassmann manifold can be used to measure the similarity between $\hat{\mathbf{R}}_i$ and $\hat{\mathbf{R}}_j$
\begin{align}
d_G(\hat{\mathbf{R}}_i,\hat{\mathbf{R}}_j)=\mathrm{VCF}_G(\mathbf{U}_{i,s},\mathbf{U}_{j,s})=\frac{\mathrm{vol}_{2s}([\mathbf{U}_{i,s},\mathbf{U}_{j,s}])}{\mathrm{vol}_s(\mathbf{U}_{i,s})\mathrm{vol}_s(\mathbf{U}_{j,s})}
\end{align}
where
\begin{align}
\text{vol}_s(\mathbf{U}_{q, s}) = \prod_{l=1}^{s} \delta_l
\end{align}
and $\delta_1, \delta_2, \dots, \delta_s \geq 0$ are the singular values of $\mathbf{U}_{q, s}$. Since $\mathbf{U}_{q, s}$ is an orthonormal basis, by substituting $\text{vol}_s(\mathbf{U}_{q, s}) = \prod_{l=1}^{s} \delta_l = 1$, the distance between two subspaces can be rewritten as
\begin{align}
d_G(\hat{\mathbf{R}}_i, \hat{\mathbf{R}}_j) = \text{VCF}_G(\mathbf{U}_{i, s}, \mathbf{U}_{j, s}) = \text{vol}_s([\mathbf{U}_{1, s}, \mathbf{U}_{2, s}])
\end{align}
Finally, the CCM can be estimated using the following optimization model
\begin{align}
\min_{\mathbf{U}_{\hat{\mathbf{R}}, s} \in G_Q(s, MN)} \sum_{q=1}^{Q} w_q \, \text{VCF}_G(\mathbf{U}_{q, s}, \mathbf{U}_{\hat{\mathbf{R}}, s})
\end{align}
where $\hat{\mathbf{R}} = \mathbf{U}_{\hat{\mathbf{R}}} \bm{\Lambda}_{\hat{\mathbf{R}}} \mathbf{U}_{\hat{\mathbf{R}}}^H$ and $\mathbf{U}_{\hat{\mathbf{R}}, s}$ is the $s$-dimensional subspace of $\mathbf{U}_{\hat{\mathbf{R}}}$. Although the VCF only depends on the eigenvector, to obtain a complete covariance matrix, it is still necessary to ensure that the eigenvalue matrices $\bm{\Lambda}_{\hat{\mathbf{R}}, s}$ and $\bm{\Lambda}_{\hat{\mathbf{R}}}$ are updated. 

\begin{algorithm}
\caption{BGVCF-STAP Algorithm.} 
\label{algorithm1} 
\begin{algorithmic}[1] 
\REQUIRE  STAP data $\mathbf{x}$, number of samples $q$, gradient descent step size $\eta$
\ENSURE  STAP suppression results $\mathbf{y}$
\vspace{0.2cm}  
\STATE Calculate THPD matrix $\mathbf{R}_q$
\STATE Establish Brauer cluster boundary limit $T_B$
\STATE Obtain clutter training samples $\hat{\mathbf{R}}_q$
\STATE Eigenvalue decomposition on $\hat{\mathbf{R}}_q = \mathbf{U}_q \bm{\Lambda}_q \mathbf{U}_q^H$
\STATE Select eigenvectors from $\mathbf{U}_q$ to form $\mathbf{U}_{q, s}$, $\bm{\Lambda}_{q, s}$ and initialize $\hat{\mathbf{R}}$
\STATE Update $\mathbf{U}_{\hat{\mathbf{R}}, s}$ using gradient descent $\mathbf{U}_{\hat{\mathbf{R}}, k} \leftarrow \mathbf{U}_{\hat{\mathbf{R}}, k} - \eta \frac{\partial \text{VCF}}{\partial \mathbf{U}_{\hat{\mathbf{R}}, k}}$
\STATE After the update, recalculate the covariance matrix $\hat{\mathbf{R}} = \mathbf{U}_{\hat{\mathbf{R}}, s} \bm{\Lambda}_{\hat{\mathbf{R}}, s} \mathbf{U}_{\hat{\mathbf{R}}, s}^H$
\STATE Space–time filter $\mathbf{w} = \frac{\hat{\mathbf{R}}^{-1} \mathbf{s}}{\mathbf{s}^H \hat{\mathbf{R}}^{-1} \mathbf{s}}$
\STATE $\mathbf{y} = \mathbf{w}^H \mathbf{x}$
\end{algorithmic}
\end{algorithm}

\section{Experimental result and analysis}
\label{sec4}
In this section, the proposed BGVCF-STAP method is evaluated and compared with several conventional approaches, including loaded sample matrix inversion (LSMI), GIP, RMT, MAP and the GVCF-STAP method without applying the BD theorem. Both simulated data and Mountain-Top measured data are utilized to demonstrate the effectiveness of proposed algorithm. 

\subsection{Simulated data}
\begin{figure}
	\centering	
        \includegraphics[width=0.5\textwidth]{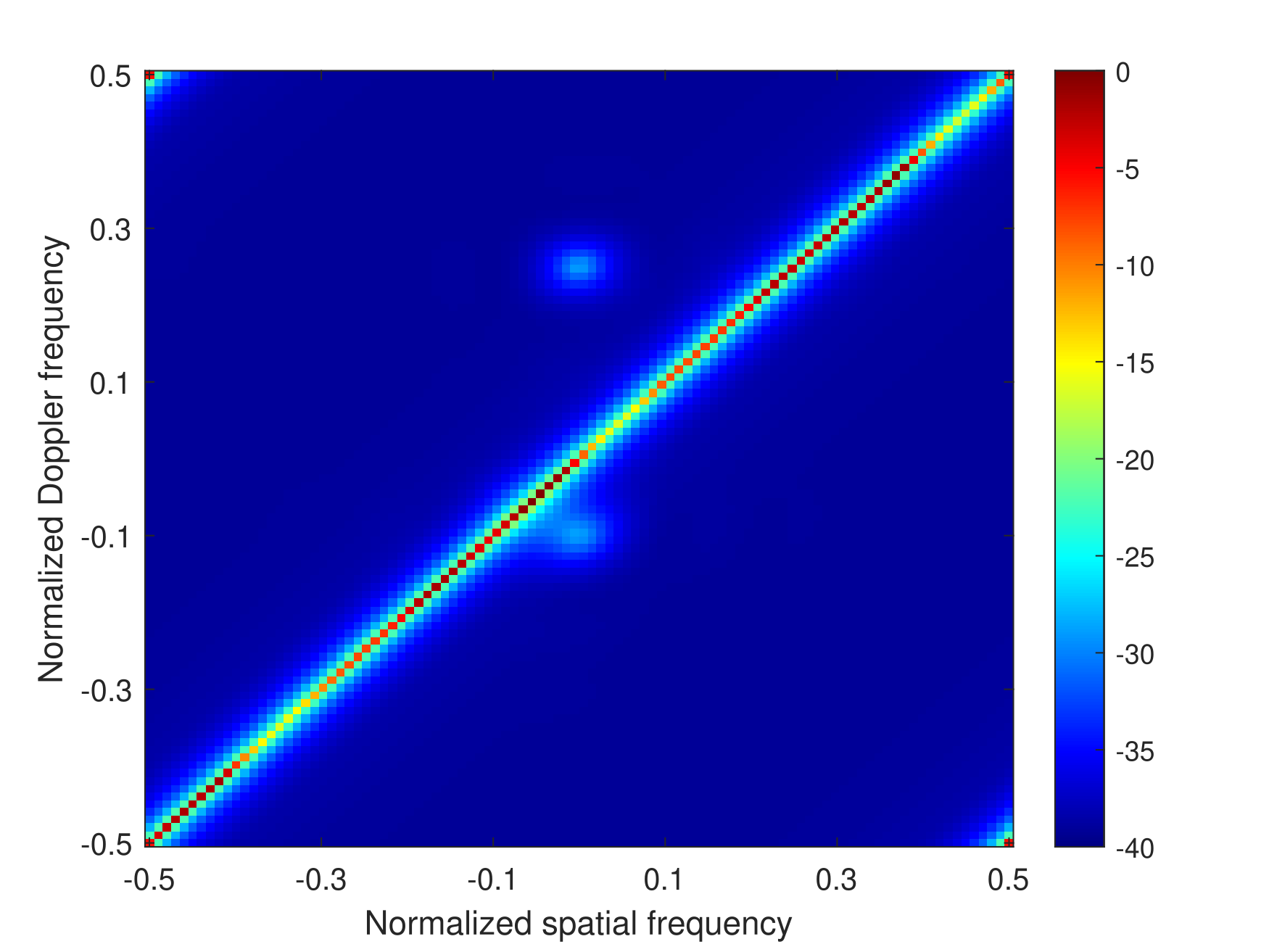}
	\caption{Spatial–temporal Capon spectrum.}
 \label{fig3}
\end{figure}

This section presents the evaluation of the proposed algorithm's performance using simulated data. The airborne radar system is modeled as a side-looking ULA, with array elements spaced at half the wavelength. The key simulation parameters are summarized in Table ~\ref{table2}. A heterogeneous clutter environment is simulated, and to account for the effect of internal clutter motion, the range-spread target 1 is introduced in range cells 30 to 32, with normalized Doppler and spatial frequencies of 0.25 and 0, respectively. Similarly, the range-spread target 2 is introduced into range cells 39–41, with normalized Doppler and spatial frequencies of -0.1 and 0, respectively. The Spatial–temporal Capon spectrum is shown in Fig. ~\ref{fig3} .
\begin{table}[h]
\small 
    \begin{center} 
        \begin{minipage}{0.7\textwidth} 
            \centering
            \caption{Simulation Parameters.} 
            \label{table2}
            \small 
            \begin{tabular}{@{}p{0.8\textwidth}cl@{}} 
                \toprule
                \text{Parameter} & \text{Value} \\ 
                \midrule
                Carrier frequency &   1.5 GHz \\ 
                PRF &   3000 Hz \\ 
                Platform velocity &   150 m/s \\ 
                Platform height &   10 km \\ 
                Number of pulses &   12 \\ 
                Number of antenna elements &   10 \\ 
                Bandwidth &   10 MHz \\  
                Number of azimuthal clutter patches &   201 \\ 
                Clutter-to-noise ratio &   50 dB \\ 
                \bottomrule
            \end{tabular}
        \end{minipage}
    \end{center}
\end{table}

\begin{figure}[htbp]
    \centering
    \begin{subfigure}[b]{0.5\textwidth}  
        \centering
        \includegraphics[width=\linewidth]{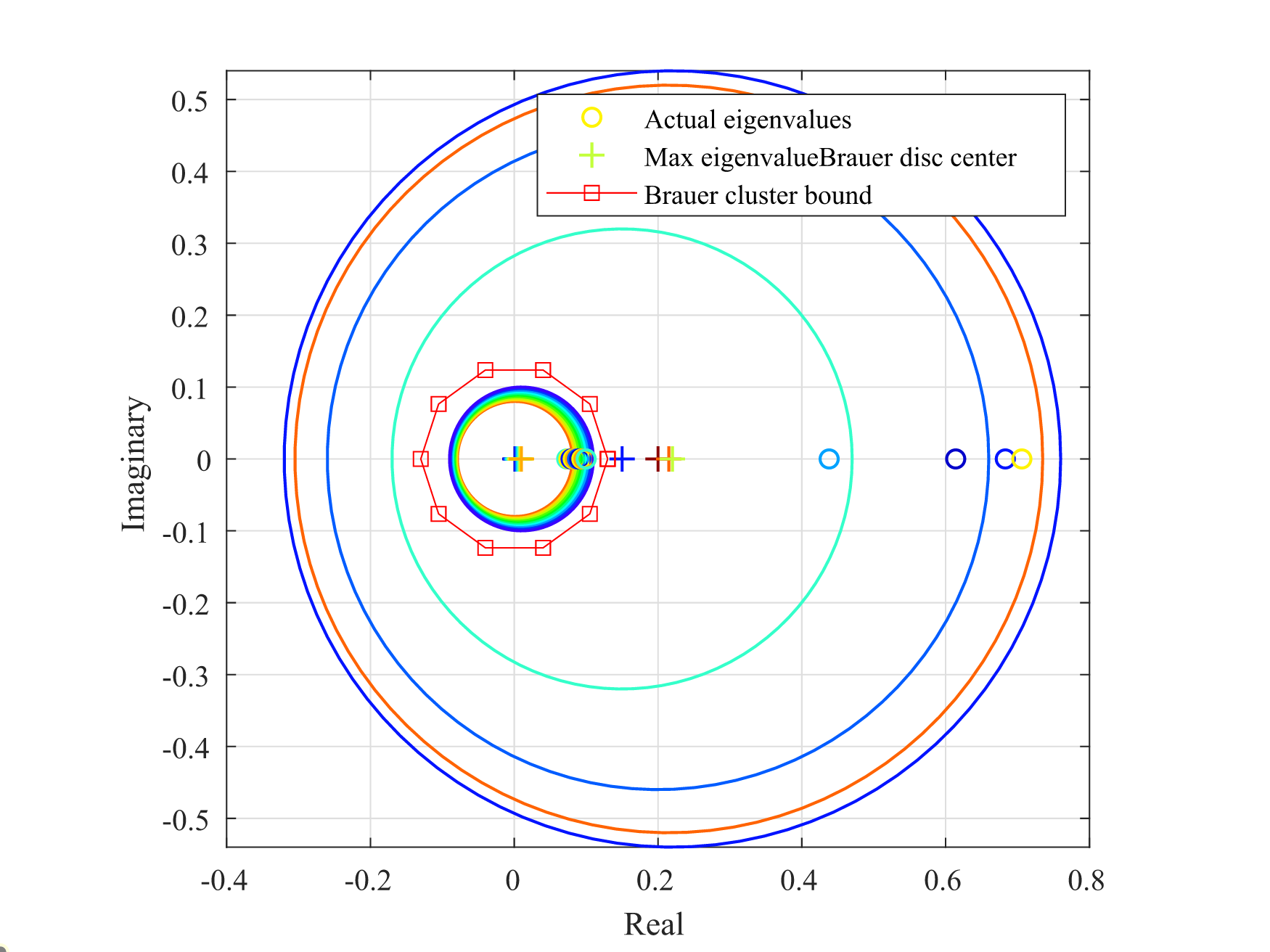}  
        \caption{}
    \end{subfigure}\\
    \begin{subfigure}[b]{0.5\textwidth}  
        \centering
        \includegraphics[width=\linewidth]{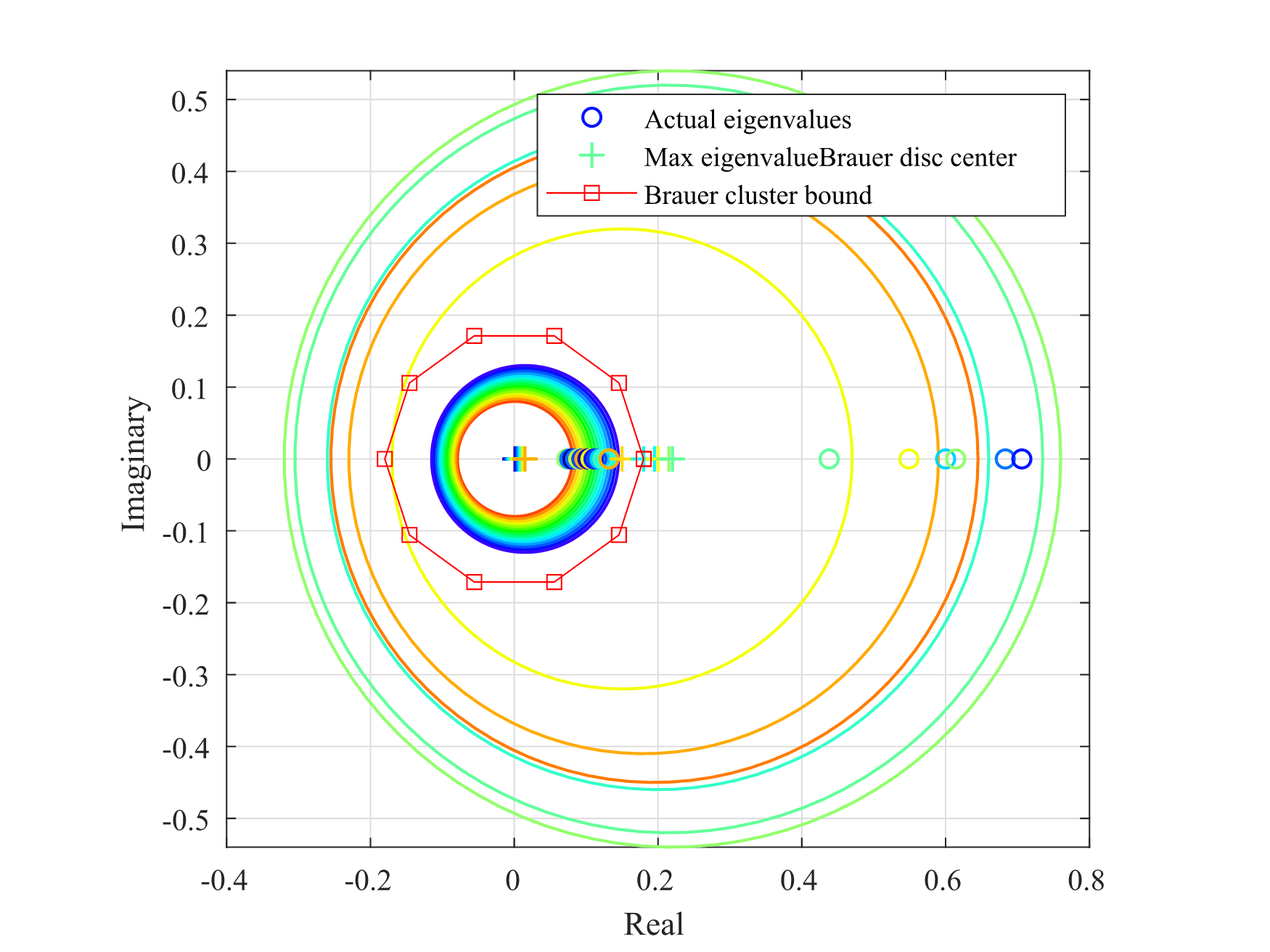}
        \caption{}
    \end{subfigure}

   \caption{Brauer bound for clutter-plus-noise and targets under different numbers of training samples. (a) number of samples = 25. (b) number of samples = 50.}
    \label{fig4}
\end{figure}

The effectiveness of the proposed techniques was assessed through two separate tests, each using 25 and 50 snapshots for analysis, respectively. The selected snapshots were transformed into THPD covariance matrices, and a clutter-plus-noise bound was determined using the BD theorem. As illustrated in Fig. ~\ref{fig4}, some points fall outside the established clutter-plus-noise bound, indicating their potential as target candidates.
In Fig. ~\ref{fig4} (a), the range cells from 15 to 39 were extracted, and four target points (30-32, 39) were detected outside the boundary. In Fig. ~\ref{fig4} (b), the range cells from 15 to 64 were extracted, and all six target points were detected.

\begin{figure}[htbp]
    \centering
    \begin{subfigure}[t]{0.4\textwidth}
        \centering
        \includegraphics[width=\textwidth]{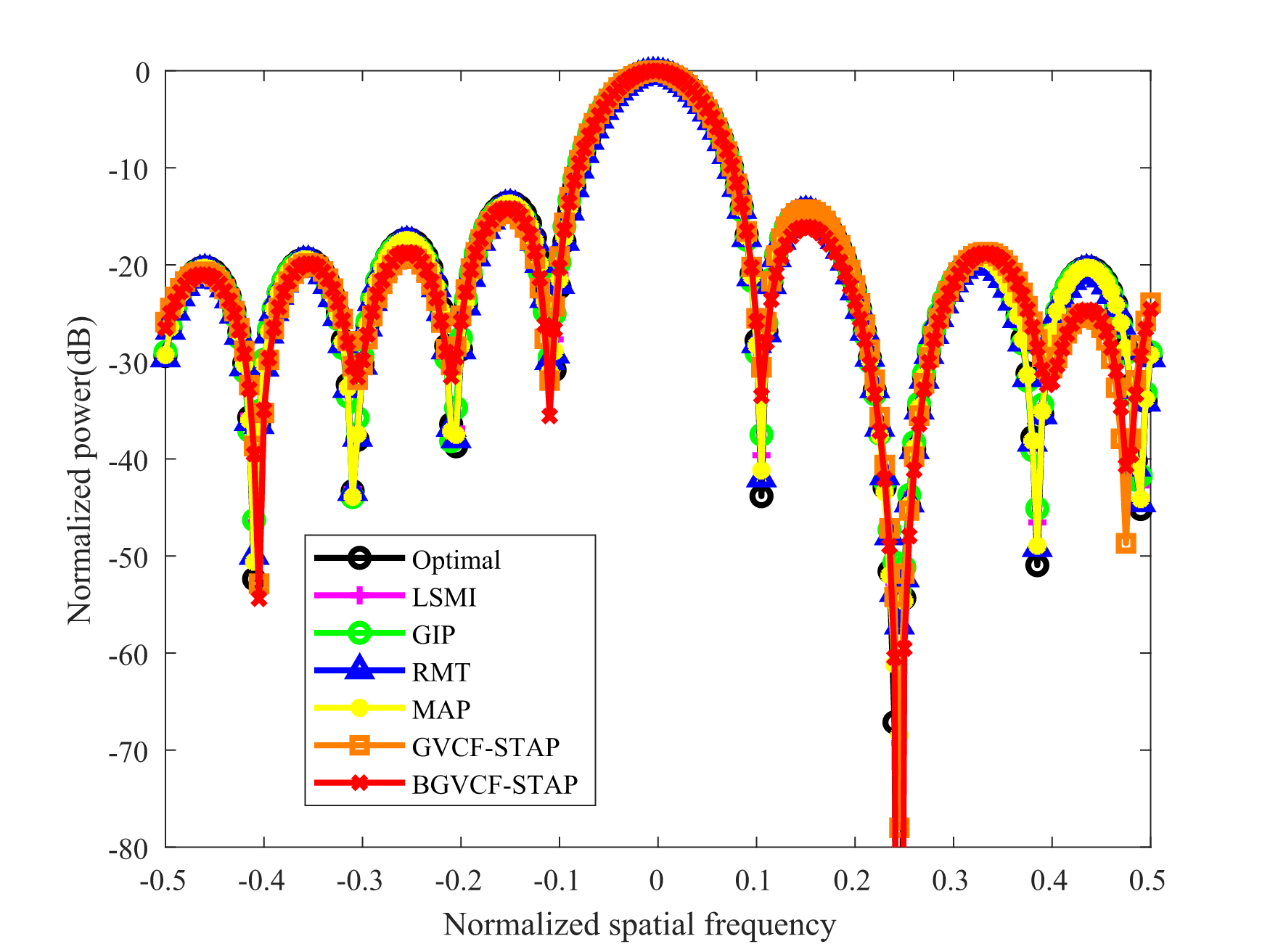}
         \caption{}
    \end{subfigure}
    \begin{subfigure}[t]{0.4\textwidth}
        \centering
        \includegraphics[width=\textwidth]{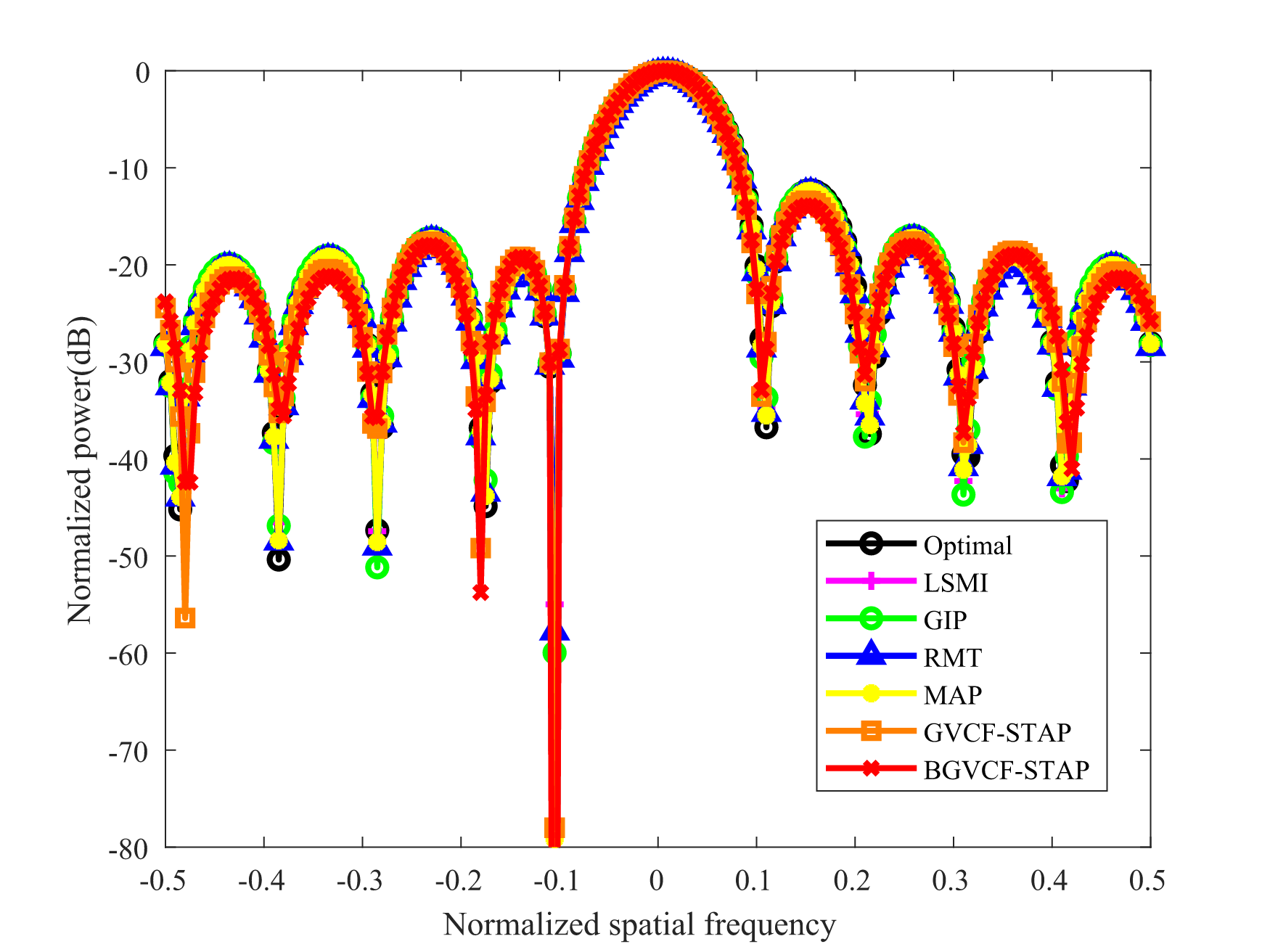}
         \caption{}
    \end{subfigure}
        \begin{subfigure}[t]{0.4\textwidth}
        \centering
        \includegraphics[width=\textwidth]{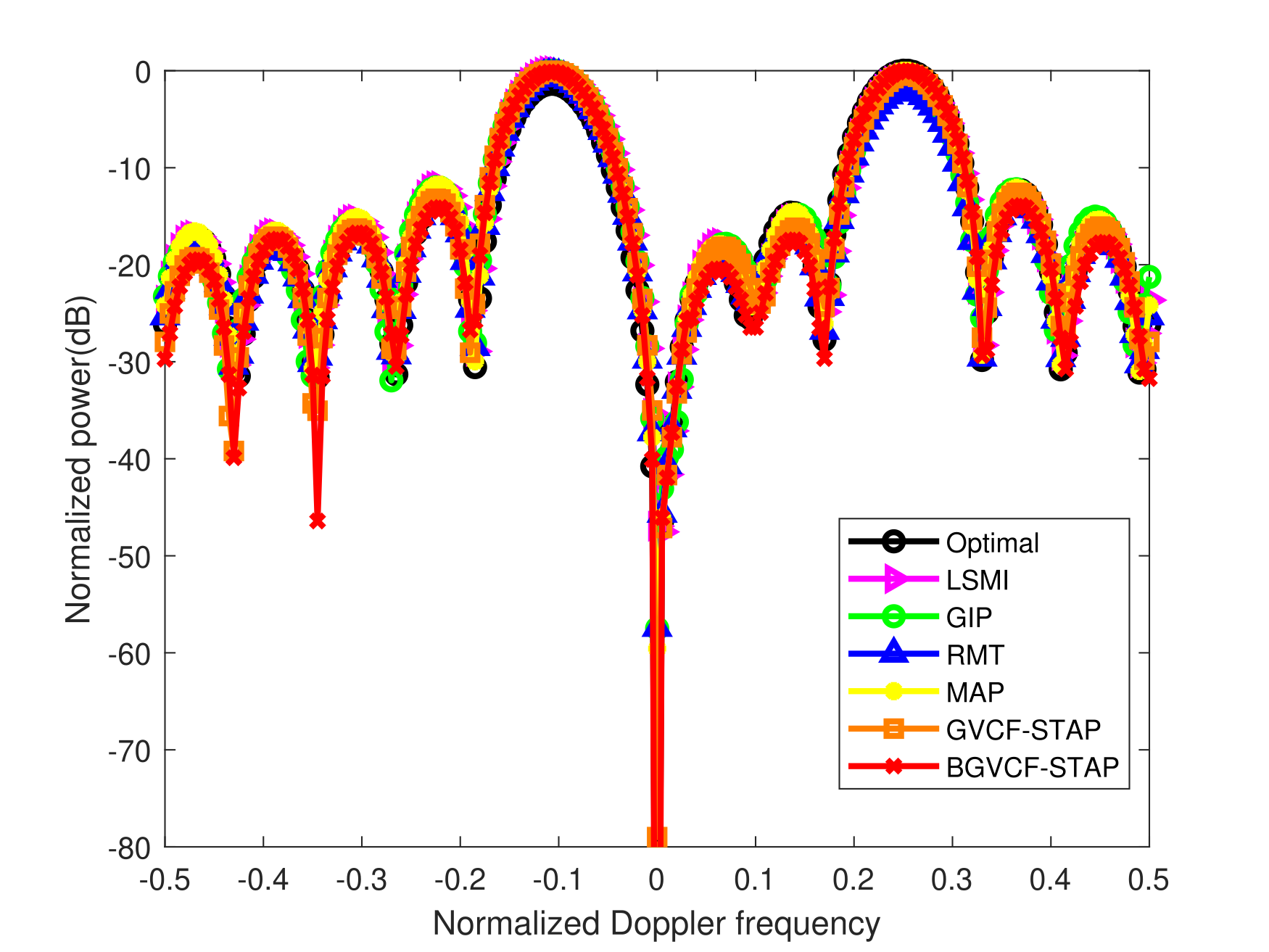}
         \caption{}
    \end{subfigure}
    \caption{Spatial and Doppler slice beampattern. (a) spatial slice beampattern of target1. (b) spatial slice beampattern of target2. (c) Doppler slice beampattern.}
    \label{fig5}
\end{figure}

Fig. ~\ref{fig5} illustrates the spatial slice and Doppler slice beampatterns, which are extracted from the space-time plane beampatterns of six algorithms at the target Doppler and spatial frequencies. These patterns are used to compare the performance details of the algorithms. Fig. ~\ref{fig5} (a) and Fig. ~\ref{fig5} (b) show the spatial slice beampatterns at the Doppler frequencies of target 1 and target 2, respectively. Traditional methods such as LSMI and GIP exhibit higher sidelobes, whereas the proposed BGVCF method achieves lower sidelobes compared to the GVCF algorithm. Additionally, the clutter notches of the BGVCF method are deeper than those of other methods. In Fig. ~\ref{fig5} (c), which depicts the Doppler slice beampattern, the proposed method also demonstrates lower sidelobes and deeper clutter notches than other algorithms. 
These results clearly indicate that the proposed method significantly outperforms others in clutter suppression, ensuring better consistency between the target location and the preset position.

\begin{figure}[htbp]
    \centering
    \begin{subfigure}[b]{0.48\textwidth}  
        \centering
        \includegraphics[width=\linewidth]{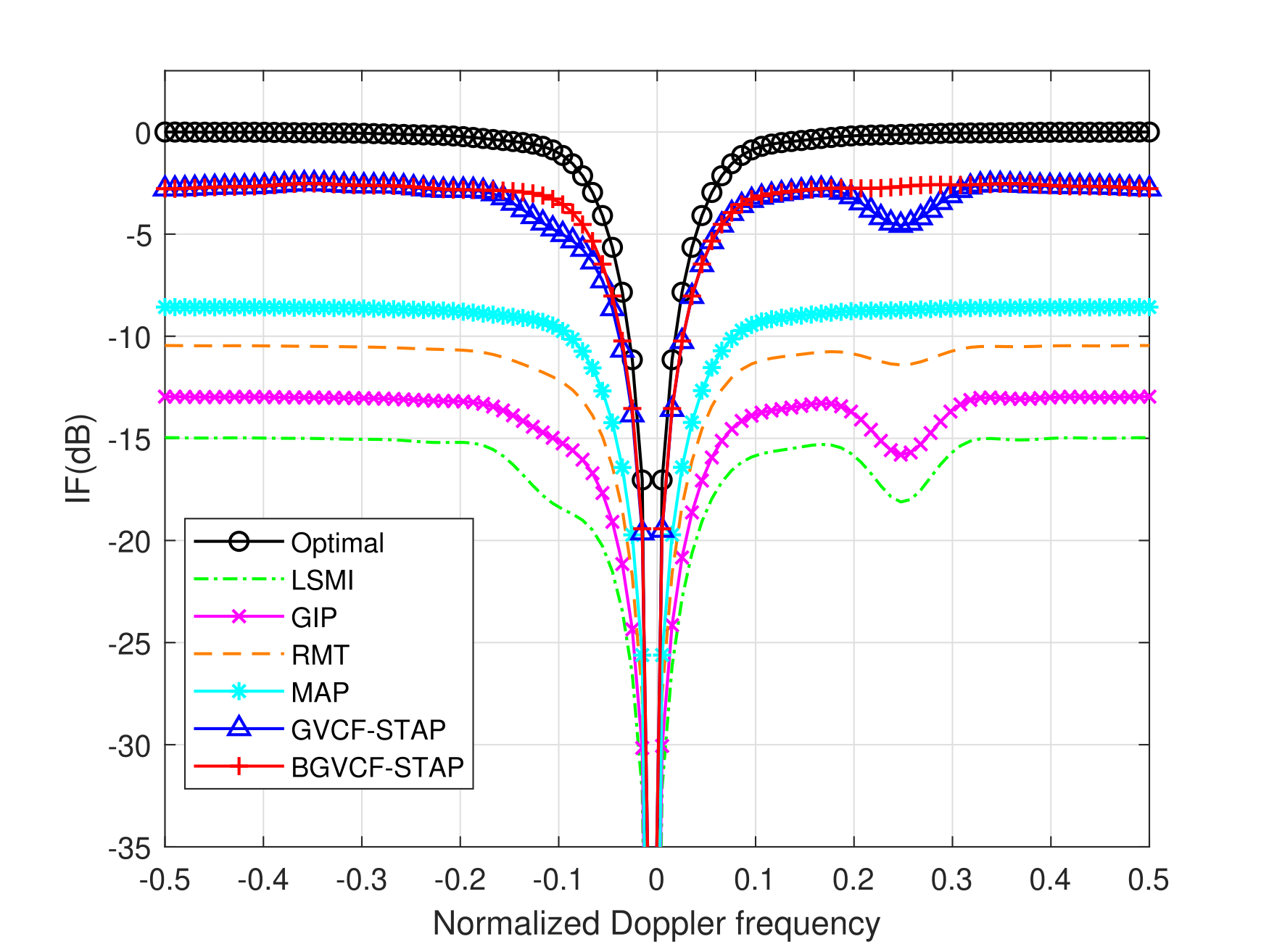}  
        \caption{}
    \end{subfigure}
    \hspace{0.01\textwidth}  
    \begin{subfigure}[b]{0.48\textwidth}  
        \centering
        \includegraphics[width=\linewidth]{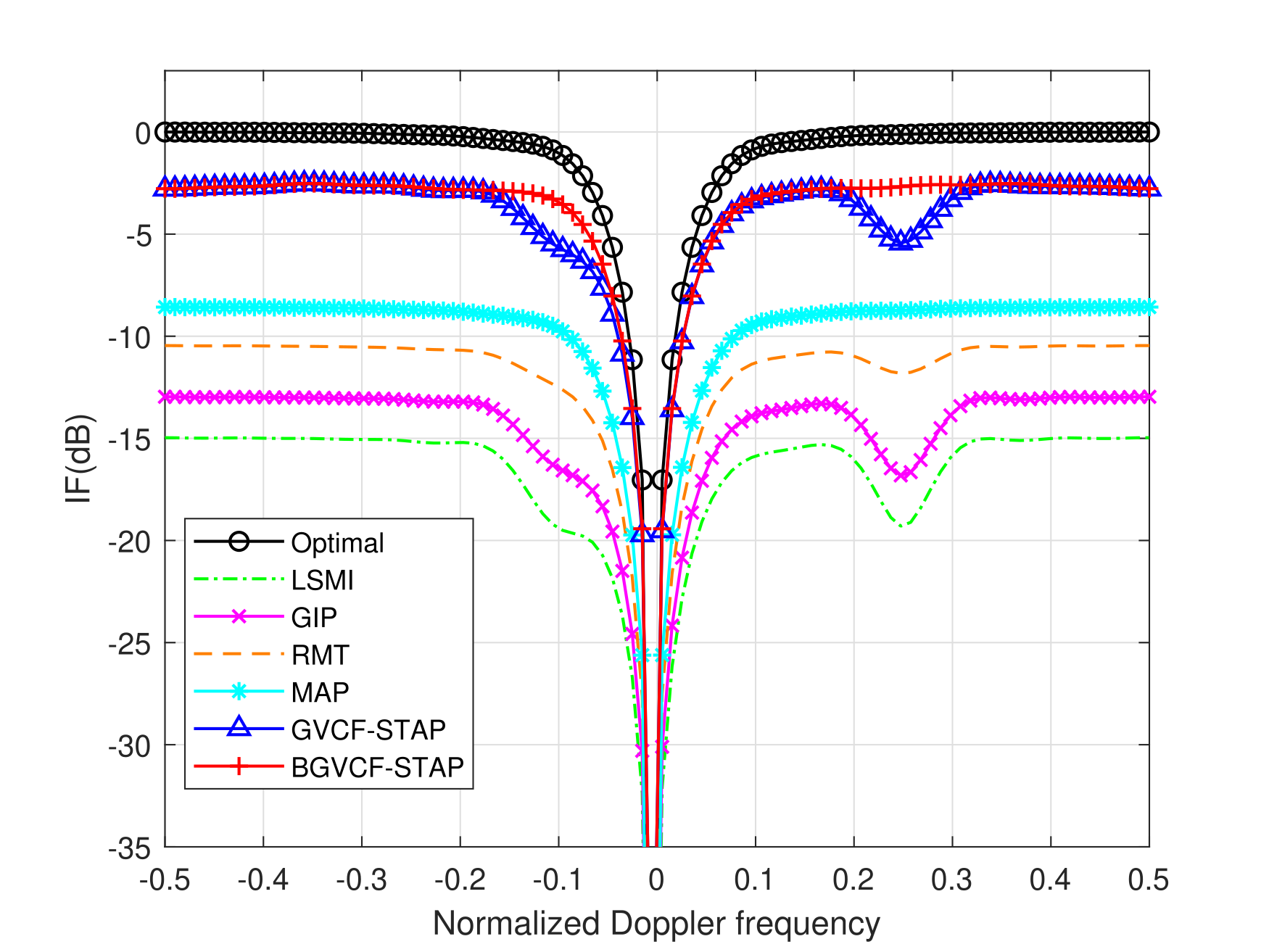}
        \caption{}
    \end{subfigure}

    \caption{Comparison of IF under different SNR. (a) SNR = 5dB. (b) SNR = 10dB.}
    \label{fig6}
\end{figure}

The improvement factor (IF) curve \cite{RN85} of different algorithms are shown in Fig. ~\ref{fig6}. The IF curve of the proposed BGVCF-STAP algorithm closely approximates the optimal curve, outperforming the MAP method by more than 5 dB and exceeding the LSMI method by over 12 dB. This is attributed to its effective utilization of the internal geometric structure of the data, without requiring any prior knowledge of the clutter distribution. GVCF-STAP algorithm suffers from performance degradation due to the presence of targets in the training sample.
 
\begin{figure}
	\centering	
        \includegraphics[width=0.5\textwidth]{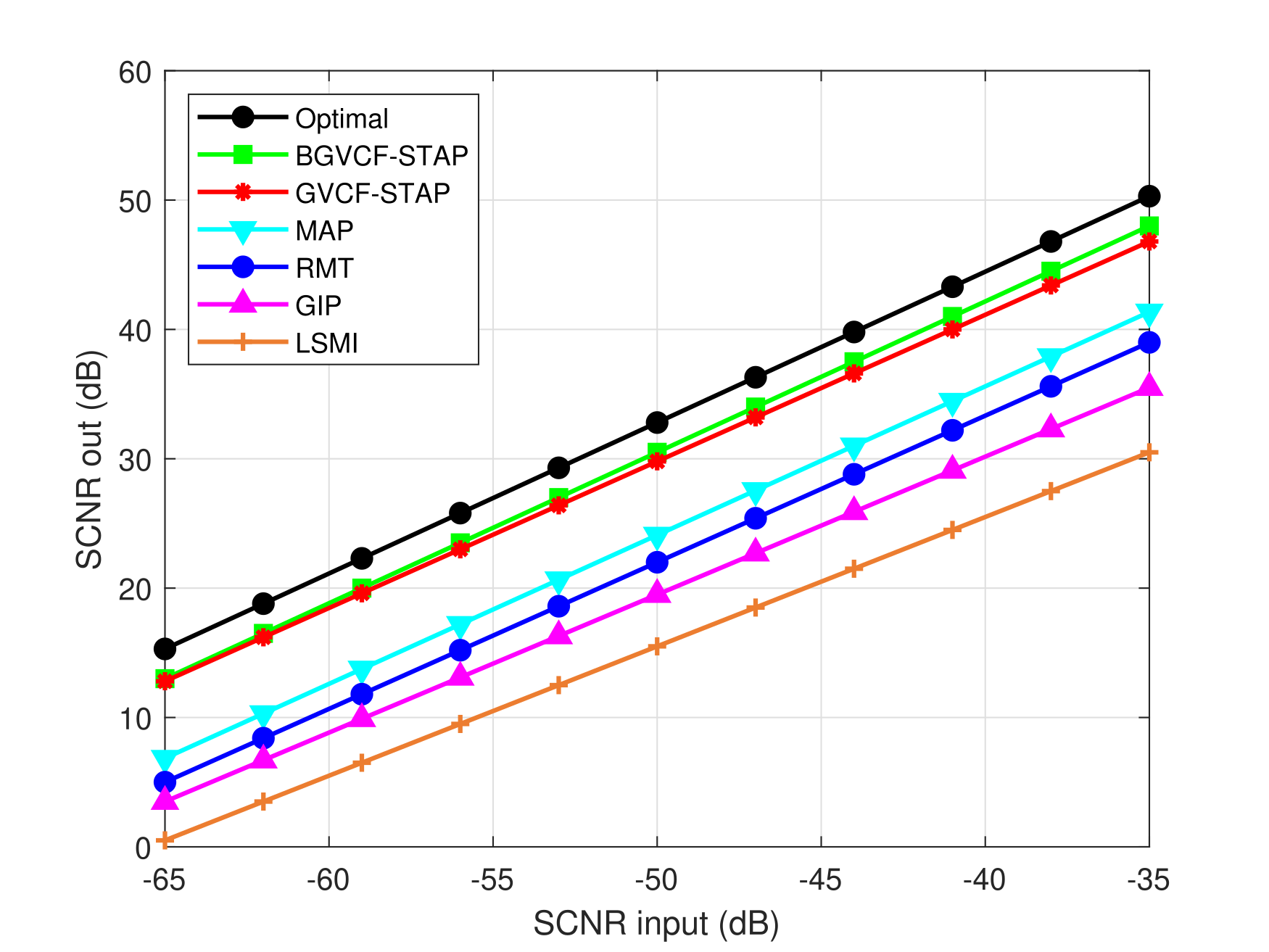}
	\caption{ Clutter suppression performance under different input SCNR.}
 \label{fig7}
\end{figure}

Fig. ~\ref{fig7} shows the output signal-to-clutter-plus-noise ratio (SCNR) \cite{RN47} of different algorithms under varying input SCNR conditions. 
BGVCF-STAP performs the best, closely approximating the optimal case. Due to the presence of targets in the training samples, GVCF-STAP performs slightly worse than BGVCF-STAP as SCNR increases. MAP ranks third, as it effectively models the actual clutter distribution. Both RMT and GIP methods are able to suppress interference from range-spread targets, but they perform worse than MAP, while the LSMI method shows the poorest performance. Overall, BGVCF-STAP and GVCF-STAP are particularly effective under low input SCNR conditions, making them suitable for scenarios that require strong clutter suppression.

\subsection{ Measured data}
\begin{figure}
	\centering
	\includegraphics[width=0.5\textwidth]{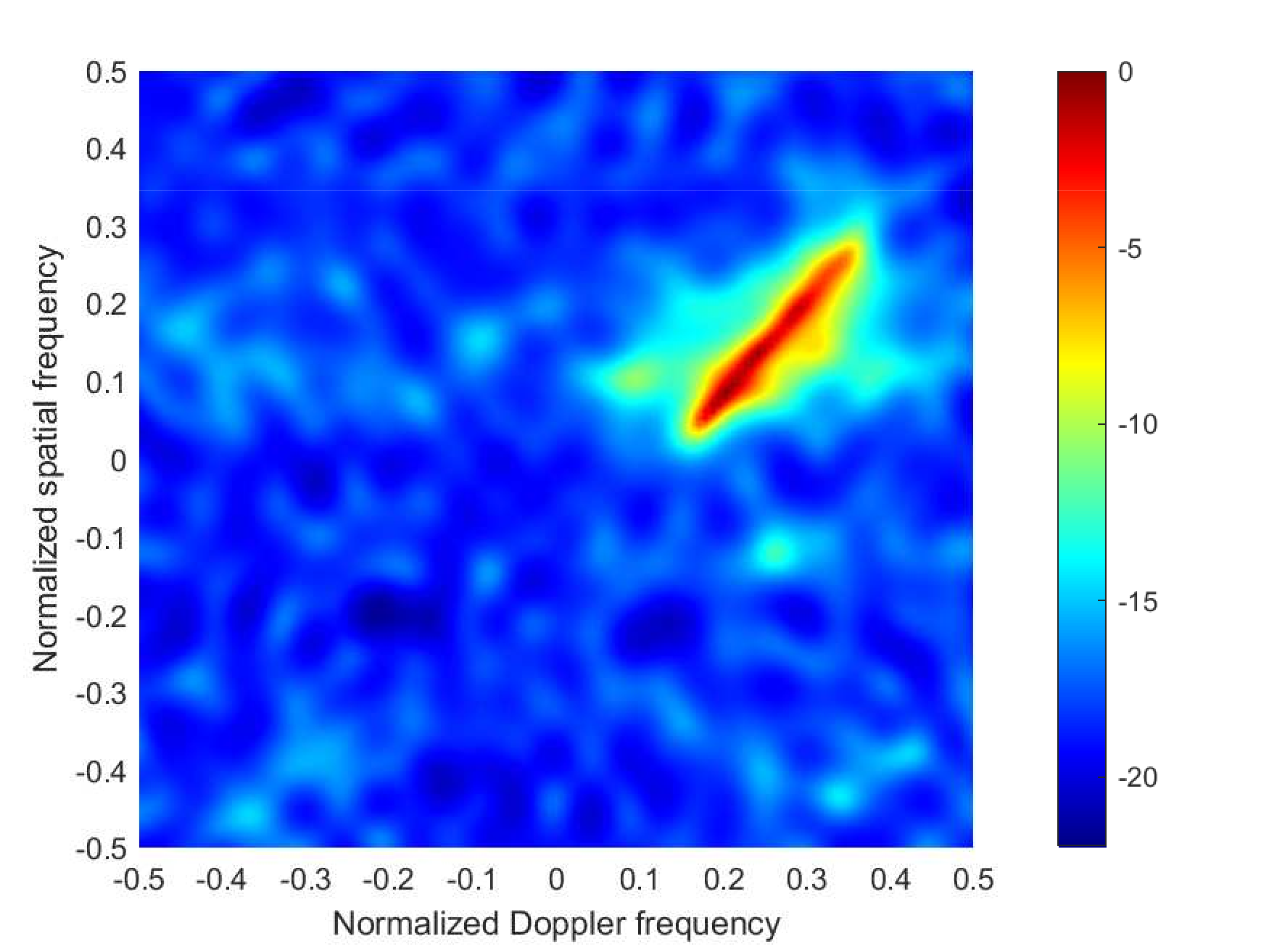}
	\caption{ Estimated clutter spectrum using all 403 range snapshots.}
 \label{fig8}
\end{figure}

We then apply the BGVCF-STAP algorithm to the publicly available Mountain-Top dataset \cite{543572}. This dataset includes a receive antenna array with 14 elements, where the echoes from 16 coherent pulses are grouped to form a CPI data set. The system used a 500 kHz linear frequency-modulated pulse for transmission, with a pulse width of 100 µs and a pulse repetition frequency (PRF) of 625 Hz. For our experiment, we utilize the data file t38pre01v1 CPI6, which contains 403 range snapshots. 
In this dataset, the target is positioned in the 147th range cell, with a normalized Doppler frequency of 0.25 and a normalized spatial angle of -15$^\circ$. The clutter power spectrum, estimated from all 403 range cells, is shown in Fig. ~\ref{fig8}.

\begin{figure}
	\centering
	\includegraphics[width=0.5\textwidth]{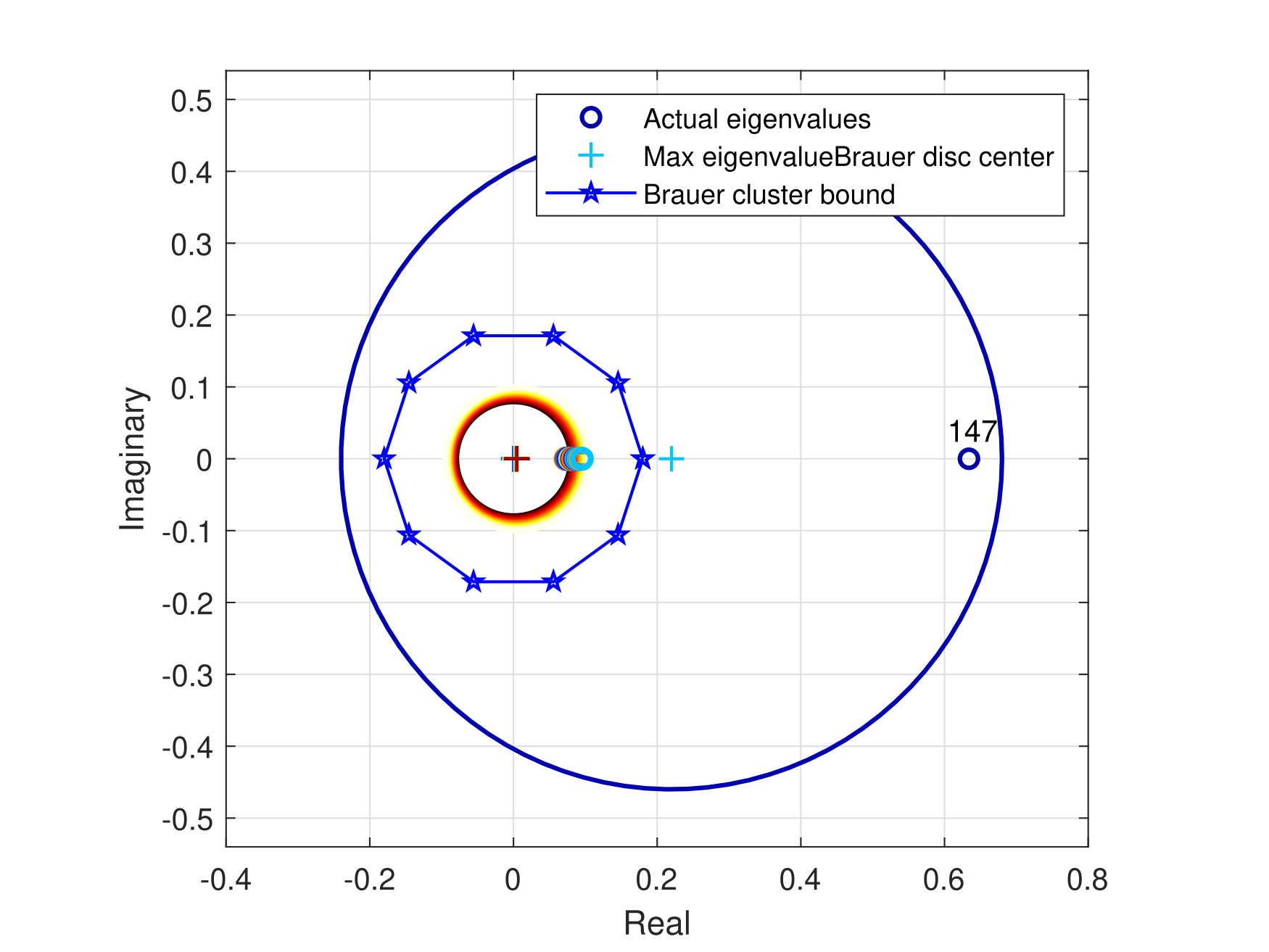}
	\caption{Brauer bound for clutter-plus-noise and target in Mountain-Top data. }
 \label{fig9}
\end{figure}

\begin{figure}
	\centering
	\includegraphics[width=0.5\textwidth]{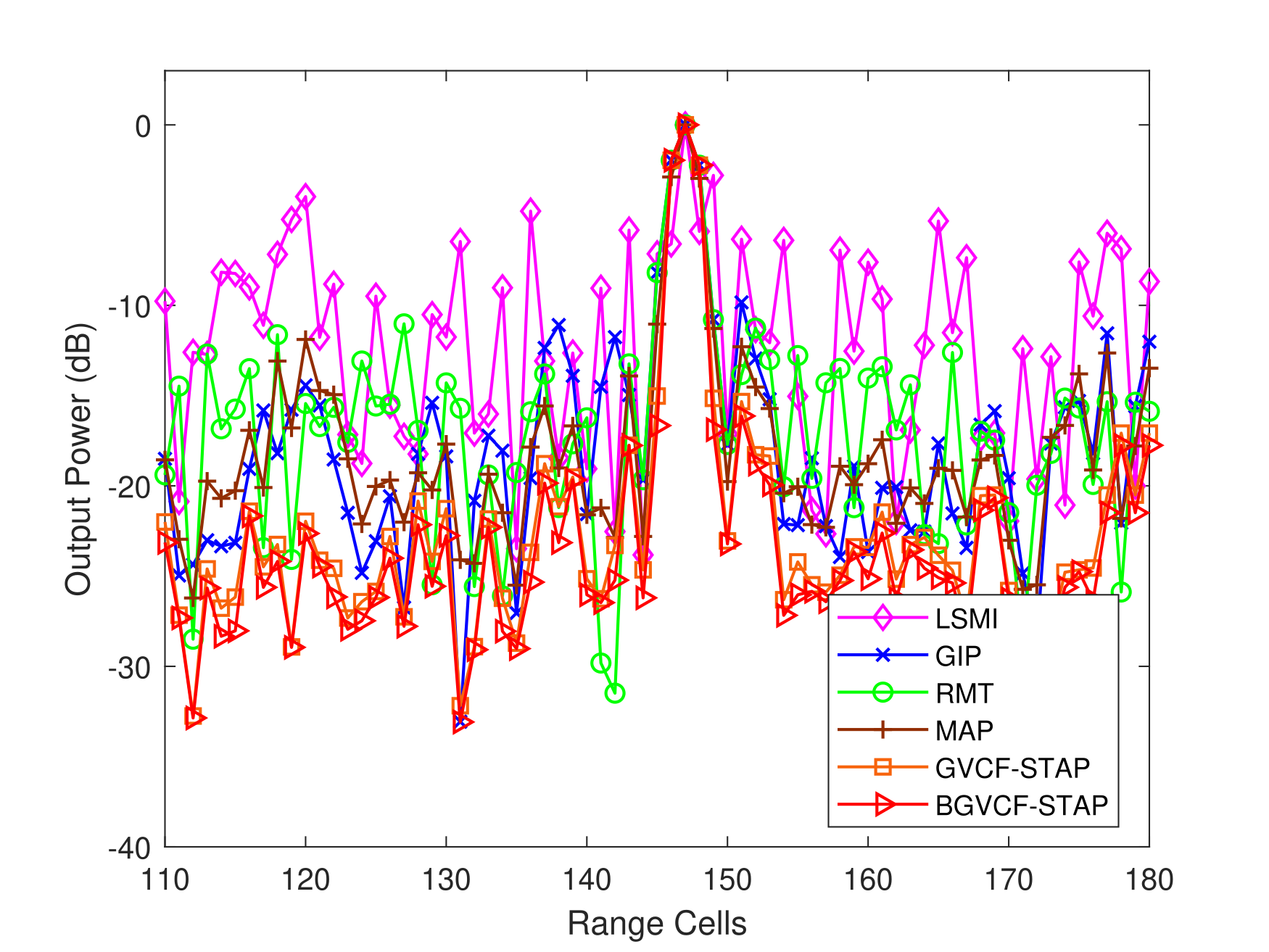}
	\caption{ Output power of different methods.}
 \label{fig10}
\end{figure}

We used the BD theorem to establish bounds for the measured data in the 110th to 180th range cells, resulting in the outcomes shown in Fig. ~\ref{fig9}. As observed, a small circle appears outside the Brauer bound for clutter-plus-noise and target, indicating that a target is detected at the 147th range cell. This finding is consistent with the known information: the target in the Mountain-Top data is located at the 147th range cell with a normalized Doppler frequency of 0.25.
In this paper, we employ a sliding window approach. For each CUT, we select 40 surrounding cells as training samples. Additionally, two guard cells are placed on both sides of the CUT. Fig. ~\ref{fig10} illustrates the STAP output power across the 110th to the 180th range cells. As depicted in the figure, the highest output power is observed at the 147th range cell, where the target is located. In contrast, the output power in the neighboring range cells, processed by the proposed BGVCF-STAP algorithm, remains low, thereby enhancing target detection.

\section{Conclusion} \label{sec5}
This paper presents a novel STAP clutter suppression method based on the BD theorem and GVCF. 
The proposed method uses a THPD covariance matrix to smooth the data, reduce variance, and enhance the consistency between the covariance matrix and the actual signal model, thereby improving estimation accuracy and algorithm performance. By applying the Brauer clustering boundary derived from the BD theorem, it effectively distinguishes range-spread targets from clutter.
The core of this method lies in transforming the clutter suppression problem into a target-clutter discrimination problem on the Grassmann manifold. By leveraging the VCF and gradient descent optimization, the method utilizes the geometric structure of the manifold space to estimate the CCM. Simulations and measured data validate the applicability of this approach in complex and heterogeneous environments. The results demonstrate that the proposed method outperforms existing techniques in clutter suppression performance.
In summary, this study introduces a novel clutter suppression technique via VCF on the manifold geometric. The algorithm not only enhances radar detection capabilities in challenging scenarios but also lays a foundation for further exploration of STAP methods on Grassmann manifolds.


\section*{Declaration of competing interest}
The authors declare that they have no known competing financial interests or personal relationships that could have appeared to influence the work reported in this paper.

\section*{Acknowledgments}
This work was supported by grants from the National Natural Science Foundation of China [No. 62171041]; and the Natural Science Foundation of Beijing Municipality [No. 4242011].


\bibliography{mybibfile}
\newpage
\end{document}